\documentclass[10pt]{book}
\usepackage{array}
\usepackage[utf8]{inputenc}
\usepackage[T1]{fontenc}
\usepackage{lmodern}
\usepackage{microtype}
\usepackage{amsmath,amsthm, amssymb}
\usepackage{enumerate}
\usepackage{calc}

\usepackage{color}
\usepackage{graphicx}

\newcommand{\N}{\mathbb{N}}
\newcommand{\e}{E}
\newcommand{\p}{P}

\newcommand{\argmin}{\mathop{\mathrm{arg\,min}}}

\newcommand{\tb}{{\frac{t}{b}}}
\newcommand{\utb}{{\frac{1-t}{b}}}

\newtheorem{definition}{Definition}
\newtheorem{remark}{Remark}

\newtheorem{lemma}{Lemma}
\newtheorem{theorem}{Theorem}

\begin{document}
\renewcommand{\baselinestretch}{1.2}
\markright{
}
\markboth{\hfill{\footnotesize\rm KARINE BERTIN AND NICOLAS KLUTCHNIKOFF
}\hfill}
{\hfill {\footnotesize\rm MINIMAX PROPERTIES OF BETA KERNEL DENSITY ESTIMATORS} \hfill}
\renewcommand{\thefootnote}{}
$\ $\par
\fontsize{10.95}{14pt plus.8pt minus .6pt}\selectfont
\vspace{0.8pc}
\centerline{\large\bf MINIMAX PROPERTIES OF BETA KERNEL DENSITY ESTIMATORS}
\vspace{.4cm}
\centerline{Karine Bertin and Nicolas Klutchnikoff}
\vspace{.4cm}
\centerline{\it  Universidad de Valpara\'{i}so, Universit\'e de Strasbourg
}
\vspace{.55cm}
\fontsize{9}{11.5pt plus.8pt minus .6pt}\selectfont

\begin{quotation}
\noindent {\it Abstract:}
In this paper, we are interested in the study of beta kernel estimators from an asymptotic minimax point of view. It is well known that beta kernel estimators are---on the contrary of classical kernel estimators---``free of boundary effect'' and
thus are very useful in practice. The goal of this paper is to prove that there is a price to pay: for very regular functions or for certain losses, these estimators are not minimax. Nevertheless they are minimax for classical regularities such as regularity of order two or less than two, supposed commonly in the practice and for some classical losses. \par

\vspace{9pt}
\par
\end{quotation}\par

\fontsize{10.95}{14pt plus.8pt minus .6pt}\selectfont
\setcounter{chapter}{1}
\setcounter{equation}{0} 
\noindent {\bf 1. Introduction}

This paper is devoted to the study of some properties concerning beta
kernel estimators.

These estimators were introduced by Chen (1999) for density estimation
with support in $[0,1]$. Indeed, contrary to classical kernel
estimators, they are (in the classical case of the estimation of a
twice differentiable density) ``free of boundary effect'': their bias
tends to $0$ (even at points $0$ and $1$) and their mean integrated
square error is of order $n^{-4/5}$.

This property has contributed to popularize their use in many applied
fields such as economy and finance. Number of papers deal with these
applications. Among others let us point out Bouezmarni and Rollin
(2003), Bouezmarni and van Bellegem (2009) and Charpentier and Oulidi
(2010).

In this article, we adopt a quite different point of view in order to
study the performance of these beta kernel estimators. We put
ourselves in an asymptotic minimax framework. In particular, our study
will not be restricted to twice differentiable density or to mean
integrated square error.


We study here the performance of beta kernel estimators in
the density model for densities belonging to H\"older spaces. We observe $n$
independent and identically distributed (i.i.d.) variables
$X_1,\ldots, X_n$ which admit the unknown density $f$ with respect to the
Lebesgue measure on $[0,1]$. Moreover, we assume that $f$ belongs
to the class of functions $\Sigma(\beta,L)$ where $\beta>0$ is a regularity
parameter and $L>0$ is a Lipshitz constant which are assumed to
be known by the statistician. The class $\Sigma(\beta,L)$ is the set of all the density functions defined on $[0,1]$ which are $m$ times differentiable such that for all $(x,y)\in[0,1]^2$:
    \[
    \left|f^{(m)}(x)-f^{(m)}(y)\right| \leq L|x-y|^{\beta-m},
    \]
where $m=\sup\{\ell\in\N: \ell<\beta\}$.

We measure the quality of the estimators of the unknown density function $f$ with a risk in $L^p$ loss. More precisely,
if $\tilde f_n$ is an arbitrary estimator, we define its risk over
$\Sigma(\beta,L)$, for $p\geq1$, by
  \[
  R_n(\tilde{f}_n,\beta,L) = \sup_{f\in\Sigma(\beta,L)}
  R_n(\tilde{f}_n,f)
  \]
  where
\[
	R_n(\tilde{f}_n,f)=\left(\e_f^n\left(\|\tilde{f}_n-f\|_p^p
\right)\right)^{\frac{1}{p}},
\]
and $\e_f^n$ is the expectation with respect to the law of $(X_1,\ldots,X_n)$.
The minimax rate of convergence on $\Sigma(\beta,L)$ is defined as
$r_n(\beta, L) = \inf_{\tilde{f}_n} R_n(\tilde{f}_n,\beta,L)$
where the infimum is taken over all the estimators. The asymptotic of $r_n(\beta,L)$ is well-known up
to a constant (see Ibragimov and Hasminskii (1981)) and is of order
$\varphi_n(\beta)=n^{-\frac{\beta}{2\beta+1}}$. We are then interested in knowing if beta kernel estimators are optimal estimators in minimax framework, in other words, if they converge at the rate $\varphi_n(\beta)$.

The advantage of this approach with respect to the classical one (second order regularity and mean integrated square error) is
that our study is quite precise concerning the understanding of the
expression ``free of boundary effect''. Our first result, Theorem~\ref{theo:N0},
illustrates that, for regularities less or equal than 2 and for $L^p$
losses with $p< 4$, it is possible to construct an optimal beta kernel
estimator. Of course, this is linked to the fact that, even on the
boundary of $[0,1]$, the bias term tends to $0$. Nevertheless, our
second result, Theorem~\ref{theo:N1}, shows that for higher regularities
($\beta>2$), even if this bias term tends  to $0$, the order of
convergence is not good.This leads to the impossibility to
construct an optimal beta kernel estimator in that case. To
conclude our study, we show, in Theorem~\ref{theo:N2}, that
even if the bias term is of the good order (for $\beta\leq 2$), for $L^p$
losses such that $p\geq 4$ the variance term is not good. This leads
again to the impossibility to construct an optimal beta
kernel estimator.

Finally, let us point out that our results clarify the conditions for
using beta kernel estimator by showing some intrinsic limitations of
these estimators. Nevertheless they have good properties of
convergence for small regularities as it is often supposed for
practical purposes. This work is a first step in the study of beta
kernel estimators from an asymptotic minimax point of view. There is still a lot to do in order to complete
this study such as finding data-driven methods to chose the bandwidth
of these estimators (cross validation methods or Lepski type
procedures).

In Section~2,  we introduce our estimators and give our main results. Section~3 is devoted to the proofs.

\par
\setcounter{chapter}{2}
\setcounter{equation}{0} 
\noindent {\bf 2. Main results}

Before stating our main results, let us recall the definition of beta kernel estimators.
\begin{definition}
  For all $b\in(0,1)$ and $t\in[0,1]$, let us introduce the following
  density:
\begin{equation}\label{eq:betadens}
  K_{t,b}(x) =
  \frac{x^\tb(1-x)^\utb}{B\left(\tb+1,\utb+1\right)}I_{[0,1]}(x)
\end{equation}
  where $B(\cdot,\cdot)$ is the standard beta function and $I$ denotes the characteristic function. Following Chen
  (1999) let us introduce the associated beta kernel
  estimator:
  \[
  \qquad\hat f_b(t) = \frac{1}{n}\sum_{k=1}^n
  K_{t,b}(X_k), \quad t\in[0,1].
  \]
\end{definition}
\begin{remark}
Let us notice that $(\hat{f}_b)_{b\in(0,1)}$ defines a one-parameter
family of estimators.
Note also that the density given by (\ref{eq:betadens}) corresponds to a beta distribution of parameters ${t}/{b}+1$ and ${(1-t)}/{b}+1$.
\end{remark}
\begin{figure}[h!]
 \includegraphics[width=\textwidth]{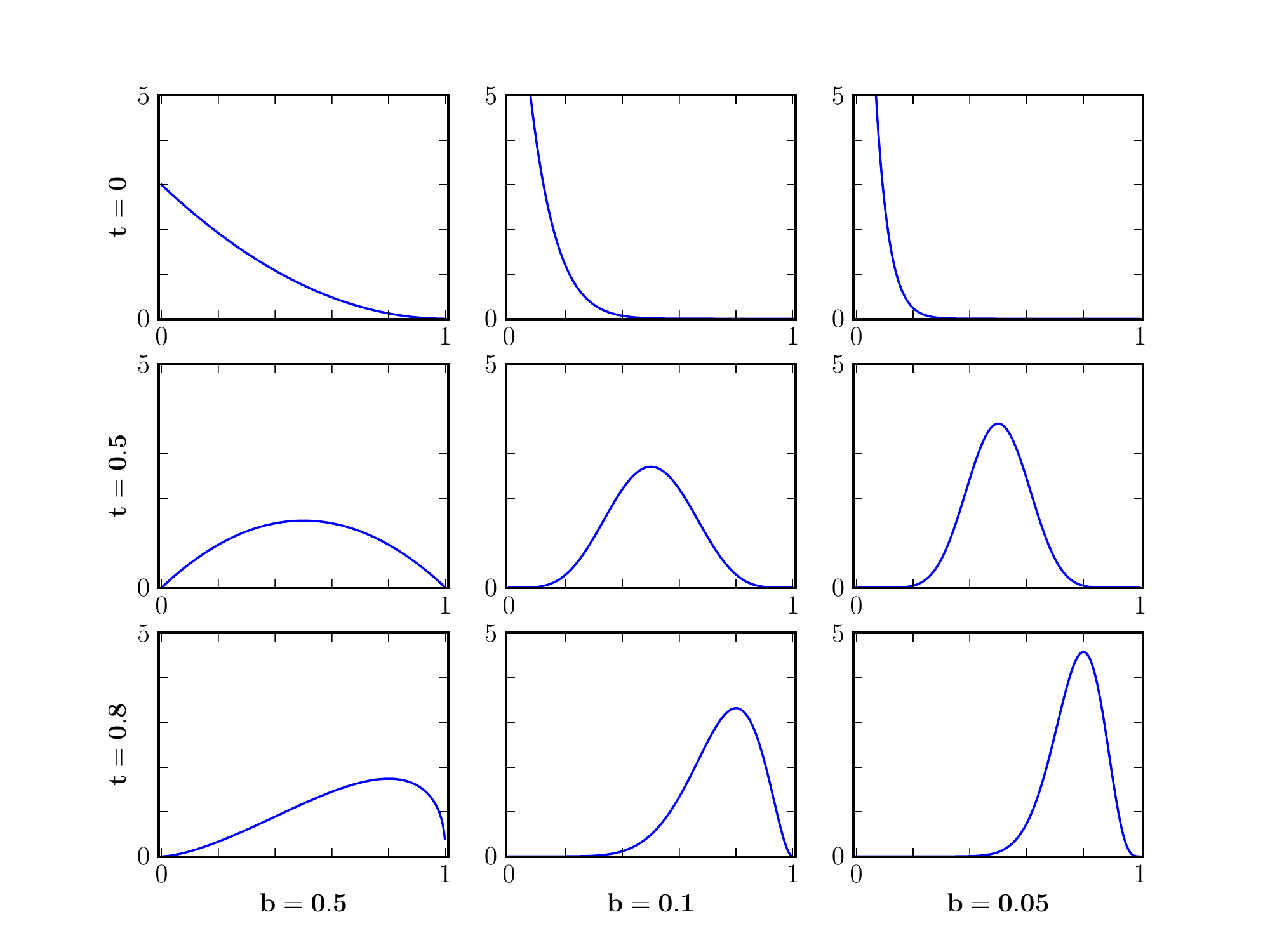}
 \caption{Some beta kernels for different values of $b$ and $t$}\label{fig:1}
\end{figure}

Figure~\ref{fig:1} represents some beta kernels $K_{t,b}$ drawn for different
values of $t$ and $b$. From left to right, $b$---which plays the role
of a bandwidth---decreases and one can observe a concentration of the
kernel in a neighborhood of $t$ which is the mode.

We will give three results about beta kernel estimators. The first one is
a positive result: for $0<\beta\leq 2$ an $1\leq p<4$ there exists
a beta kernel estimator with properly chosen \emph{bandwidth} which
achieves the minimax rate of convergence on $\Sigma(\beta,L)$ in
$L^p$-loss.


\begin{theorem}
  \label{theo:N0}
  Set $1\leq p < 4$ and $0<\beta\leq2$. Set $b_n=cn^{-\frac{2}{2\beta+1}}$,
  where $c$ is a positive constant. Then the estimator $\hat{f}_{b_n}$ achieves the rate $r_n(\beta,L)$.
  More precisely:
  \[
  \limsup_{n\to\infty}
  \frac{R_n(\hat{f}_{b_n},\beta,L)}{r_n(\beta,L)}<\infty.
  \]
\end{theorem}

\par
The following theorem shows that, for regularities larger than $2$ and most of all losses, beta kernel density estimators are not minimax.
\begin{theorem}
  \label{theo:N1}
  Set $p\ge 2$ and $\beta>2$.
  Then the family of estimators $\hat{f}_{b_n}$ satisfies, for all sequence $(b_n)_{n\in\mathbb{N}}$ in $(0,1)$,
  \[
  \liminf_{n\to\infty}
  \frac{R_n(\hat{f}_{b_n},\beta,L)}{r_n(\beta,L)}=\infty.
  \]
\end{theorem}

The next theorem proves, for regularities less or equal than 2,  contrary to what happens in Theorem~\ref{theo:N0}, that the beta kernel density estimator is not minimax if the risk is measured with $L_p$ losses with $p\geq 4$.

\begin{theorem}
  \label{theo:N2}
  Set $p\ge 4$ and $0<\beta\leq 2$.
  Then the family of estimators $\hat{f}_{b_n}$ satisfies, for all sequence $(b_n)_{n\in\mathbb{N}}$ in $(0,1)$,
  \[
  \liminf_{n\to\infty}
  \frac{R_n(\hat{f}_{b_n},\beta,L)}{r_n(\beta,L)}=\infty.
  \]
\end{theorem}

In the framework of Theorem~\ref{theo:N0}, the bias term of beta kernel estimators is of order $b^{\beta/2}$ and variance term $(nb^{1/2})^{-1/2}$. These orders are the same as for classical kernel estimators with $h=b^{1/2}$.
A trade off between bias and variance leads to choose $b$ of the form $cn^{-\frac{2}{2\beta+1}}$ and the associated beta kernel estimator converges to $f$ at rate $\varphi_n(\beta)$.
Let us remark that this theorem can be viewed as a generalization of the ``classical case'' which corresponds to $\beta=2$ and $p=2$.

When $\beta>2$ and $p\geq 2$, the impossibility to construct an optimal beta kernel estimator is linked with the construction of two functions which belong to $\Sigma(\beta,L)$. The first one has its bias term that cannot go to $0$ faster than the rate $b$ while the variance term of the second one is lower bounded by $(nb^{1/2})^{-1/2}$. This implies that beta kernel estimators cannot converge to $f$ at a rate faster than $\varphi_n(2)$ in this framework.

In the case of Theorem~\ref{theo:N2}, the non-optimality of beta kernel estimators is linked with the construction of a function, which belongs to $\Sigma(\beta,L)$, such its variance term is not of good order (by at least an extra $|\log b|$ factor). We exhibit a second function in $\Sigma(\beta,L)$ with a bias term of order $b^{\beta/2}$ that allows us to conclude that beta kernel estimators cannot converge at the minimax rate $\varphi_n(\beta)$.

Proofs of these three theorems are given in Section 3.

\par
\setcounter{chapter}{3}
\setcounter{equation}{0} %
\noindent {\bf 3. Proofs}
\par
In all the proofs, $C$ denotes a positive constant that can change of values from line to line.
In the following, $(a_n)\asymp (b_n)$ stands for $0<\liminf_{n\to\infty}
a_n/b_n \leq \limsup_{n\to\infty} a_n/b_n < +\infty$.

\par
{\bf 3.1 Proof of Theorem~\ref{theo:N0}}
\par
Firstly, set $1\leq p<4$, $0<\beta\leq2$, $L>0$ and $f\in\Sigma(\beta,L)$. As the function $f$ (and the number of observations $n$) is always fixed in this proof, we will denote, for simplicity, $\e$ instead of $\e_f^n$. We have:
\begin{eqnarray*}
  R_n\left(\hat f_b, f\right) & = & \left(\e\left(\|\hat f_b-f\|_p^p\right)\right)^{1/p} \\
  & = & \left(\e\left(\int_0^1\left|\hat f_b(t)-f(t)\right|^pdt\right)\right)^{1/p} \\
  & \leq & C\left(\int_0^1|B_t|^pdt+\e\left(\int_0^1|Z_t|^pdt\right)\right)^{1/p},
\end{eqnarray*}
where $B_t = \e\left(\hat f_b(t)\right)-f(t)$ and $Z_t = \hat
f_b(t) - \e\left(\hat f_b(t)\right)$.
Thus, as $p\geq 1$, we obtain:
\begin{equation}\label{eq:majorationRisk}
  R_n\left(\hat f_b, f\right) \leq
  C\left\{\left( \int_0^1|B_t|^pdt\right)^{\frac{1}{p}}+\left(\int_0^1 \e\left(|Z_t|^p\right)dt\right)^\frac{1}{p} \right\}.
\end{equation}
The proof of our theorem will be derived from two lemmas. The
first one is used to control the integrated bias term and the
second one to control the integrated moment of the centered stochastic
term $Z_t$.
\begin{lemma}
  Set $p$ and $\beta$ as in Theorem~\ref{theo:N0}. Then, for all $0<b<1$:
  \begin{equation}\label{eq:majbiais}
    \left(\int_0^1|B_t|^pdt\right)^{\frac{1}{p}} \leq C L b^{\frac{\beta}{2}}.
  \end{equation}
\end{lemma}
\begin{lemma}
  Set $p$ and $\beta$ as in Theorem~\ref{theo:N0}. Then, for all $0<b<1$:
  \begin{equation*}
    \left(\int_0^1\e\left(|Z_t|^p\right)dt\right)^{\frac{1}{p}} \leq C \left(\frac{1}{nb^{\frac{1}{2}}}\right)^{\frac{1}{2}}.
  \end{equation*}
\end{lemma}
Proofs of these lemmas are in Subsection 3.3 and 3.4. Let us complete
the proof of Theorem~\ref{theo:N0}. From Equation~(\ref{eq:majorationRisk})
and Lemmas~1 and~2 we derive:
\begin{equation*}
  R_n\left(\hat f_b, f\right) \leq
  C \left\{L b^{\frac{\beta}{2}}+\left(\frac{1}{nb^{\frac{1}{2}}}\right)^{\frac{1}{2}} \right\}.
\end{equation*}
As this equation is valid for all $0<b<1$ we deduce that, if $b_n$
denotes the quantity
$cn^{-\frac{2}{2\beta+1}}\asymp\argmin_{0<b<1}(L
b^{\beta/2}+(nb^{1/2})^{-(1/2)})$ with $c$ a positive constant, we have :
\[
R_n\left(\hat f_{b_n}, f\right) \leq C \left\{L
  b_n^{\frac{\beta}{2}}+\left(\frac{1}{nb_n^{\frac{1}{2}}}\right)^{\frac{1}{2}}
\right\}\leq Cr_n(\beta,L),
\]
that gives the result of the theorem.

\par 
{\bf 3.2 Proofs of Theorem~\ref{theo:N1} and \ref{theo:N2}}
\par
\begin{figure}[htp!]
\centering
 \includegraphics[width=\textwidth]{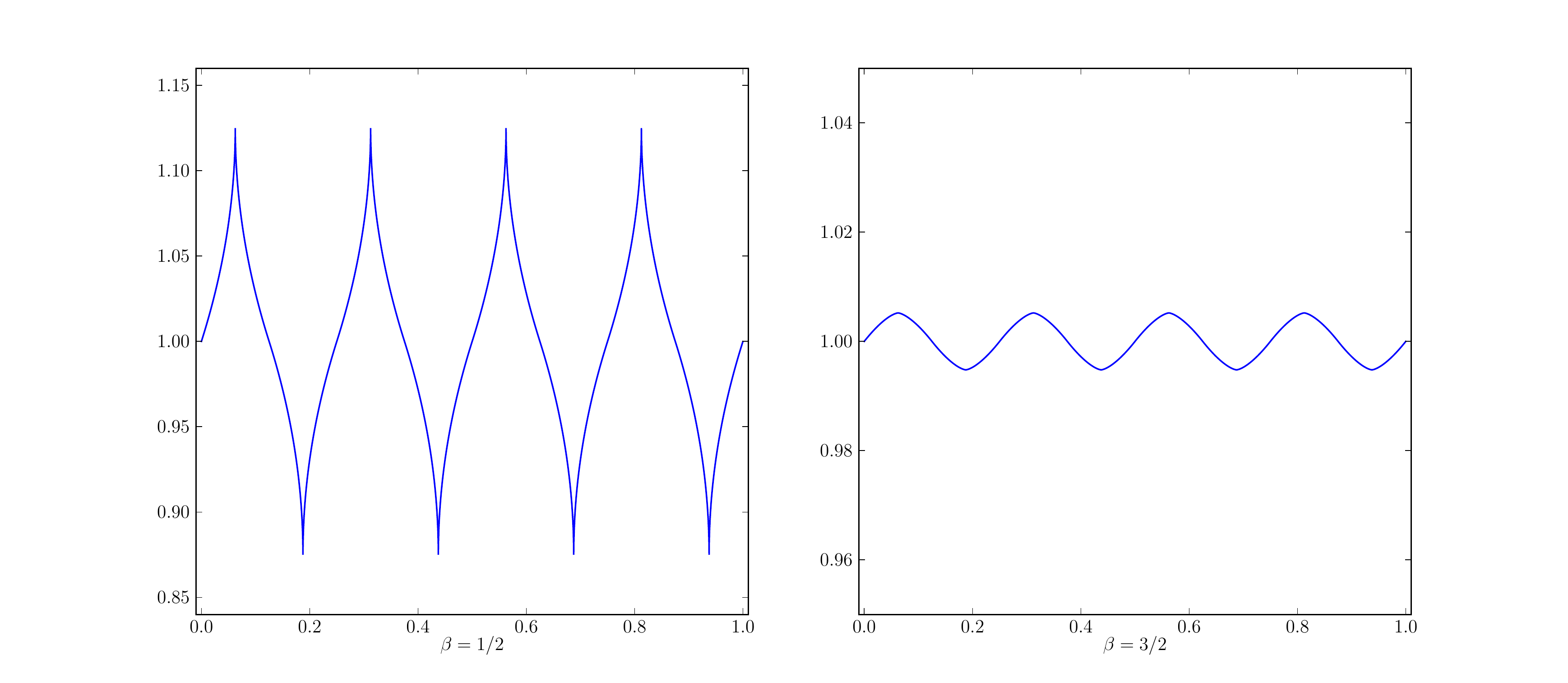}
 \caption{Function $f_\beta$: $\beta=1/2$ in the left side, $\beta=3/2$ in the right side. In both cases, $L=1$ and $b=0.0005$.}\label{fig:f1}
\end{figure}

In all these proofs we will use auxiliary functions. Let us define the first one:
\begin{equation*}
 f_0(x) = I_{[0,1]}(x).
\end{equation*}
Next, for $0<\beta\leq 2$, we define
\begin{equation*}
	f_\beta(x) = 1+L_\beta\sum_{k=1}^{2N} (-1)^{k+1} \left(\frac{1}{(4N)^\beta}-\left|x-\frac{2k-1}{4N}\right|^\beta\right)I_{[\frac{k-1}{2N},\frac{k}{2N}]}(x),
\end{equation*}
where $N$ is the integer part of $b^{-1/2}/20$ and $L_\beta=\frac{L}{2}\min(1, 1/\beta)$.  Figure~\ref{fig:f1} shows the function $f_\beta$ for two values of $\beta$. Finally we consider:
\begin{equation*}
 f_3(x)=2xI_{[0,1]}(x).
\end{equation*}

For simplicity, we will use the notation $\e_\alpha$ instead of $\e_{f_\alpha}^n$. These functions satisfy the following lemma:
\begin{lemma}\label{LEM:10}
	For all $p\geq 2$ and $b\in(0,1)$, we have:
	\begin{enumerate}[i)]
		\item The function $f_0$ belongs to $\Sigma(\beta,L)$ for all $\beta>0$ and $L>0$ and satisfies:
		\begin{equation}\label{eq:minvar_0}
		\e_0\left(\|\hat f_b-f_0\|_p^{p}\right) \geq Cd_n(b,p)\left(nb^{1/2}\right)^{-p/2},
		\end{equation}
		where
		$$d_n(b,p)=\int_{\frac{1}{2}}^{1-b}\left(\frac{1}{1-t}\right)^{\frac{p}{4}}dt.$$
		
		\item Set $0<\beta\leq 2$ and $L>0$. The function $f_\beta$ (which also depends on $b$ and $L$) belongs to $\Sigma(\beta,L)$ and satisfies:
		\begin{equation}\label{eq:minbiais_1}
		\e_\beta\left(\|\hat f_b-f_\beta\|_p^{p}\right) \geq Cb^{p\beta/2},
	\end{equation}
		
		\item The function $f_3$ belongs to $\Sigma(\beta,L)$ for all $\beta>2$ and $L>0$ and satisfies:
		\begin{equation}\label{eq:minbiais_3}
		\e_3\left(\|\hat f_b-f_3\|_p^{p}\right) \geq Cb^{p}
		\end{equation}
	\end{enumerate}
\end{lemma}

Let $\beta>2$, $L>0$ and $p\geq 2$. We have for all $b\in(0,1)$, using~(\ref{eq:minvar_0}) and~(\ref{eq:minbiais_3}):
\begin{eqnarray*}
	\sup_{f\in\Sigma(\beta,L)} R_n(\hat f_b, f) & \geq & \frac{1}{2} \left(R_n(\hat f_b, f_0)+R_n(\hat f_b, f_3)\right) \\
	& \geq & C\inf_{b\in(0,1)}\left(b+\left(\frac{1}{nb^{1/2}}\right)^{1/2}\right) \asymp \varphi_n(2).
\end{eqnarray*}
Since $\varphi_n(\beta)=o\left(\varphi_n(2)\right)$, this implies Theorem~\ref{theo:N1}.

Let $0<\beta\leq 2$, $L>0$ and $p\geq 4$. We have for all $b\in(0,1)$:
\begin{equation*}
	\sup_{f\in\Sigma(\beta,L)} R_n(\hat f_b, f) \geq \frac{1}{2} \left(R_n(\hat f_b, f_0)+R_n(\hat f_b, f_\beta)\right).
\end{equation*}
Then we have, using~(\ref{eq:minvar_0}) and~(\ref{eq:minbiais_1}):
\begin{eqnarray*}
	\sup_{f\in\Sigma(\beta,L)} R_n(\hat f_b, f) & \geq & C\left(b^{\beta/2}+\frac{d_n(b,p)}{(n b^{1/2})^{1/2}}\right)\\
	& \geq & C\left(b^{\beta/2}+\frac{|\log b|}{(n b^{1/2})^{1/2}}\right).
\end{eqnarray*}

The function $f(b)=b^{\beta/2} +|\log b|/(n b^{1/2})^{-1/2}$ attains its minimum $b^0_n$ on $(0,1)$ for $n$ large enough. The minimum satisfies $b^0_n=c_1n^{-2/(2\beta+1)}(c_2+|\log (b^0_n)|)^{4/(2\beta+1)}$ with $c_1$ and $c_2$ positive constants. It can be easily proved that the sequence $b^0_n$ tends to $0$ as $n$ tends to $\infty$. This implies that
$$
\inf_{b\in (0,1)}\left\{b^{\beta/2} + \frac{|\log b|}{(n b^{1/2})^{1/2}}\right\}=(b^0_n)^{\beta/2} + \frac{|\log b^0_n|}{(n (b^0_n)^{1/2})^{1/2}}\asymp \varphi_n(\beta) |\log (b^0_n)|^{2\beta/(2\beta+1)}.
$$
This last result implies Theorem 3.

\par 
{\bf 3.3 Proof of lemma 1}
\par
To prove inequality (\ref{eq:majbiais}), it is sufficient to prove that $|B_t| \leq C Lb^{\frac{\beta}{2}}$, for all $t\in[0,1]$.
In order to prove this last inequality, we will use the following
result (see Johnson, Kotz and Balakrishnan, 1994).
\begin{lemma}\label{LEM:JKB}
  If $\xi$ is a variable with density $K_{t,b}$ ($\xi\sim K_{t,b}$) then there exist two functions $\Delta_1$ and
  $\Delta_2$, and a constant $M>0$ such that:
  \begin{equation}
    \label{eq:20}
    \e\left(\xi\right)-t = b(1-2t) + \Delta_1(t,b)
  \end{equation}
  and
  \begin{equation}
    \label{eq:21}
    \sigma^{2}(\xi) = bt(1-t)+\Delta_2(t,b),
  \end{equation}
  where
  \begin{equation}
    \label{eq:23}
    \sup_{0\leq t\leq 1}|\Delta_j(t,b)| \leq Mb^2 \quad\text{if}\quad j=1,2.
  \end{equation}
\end{lemma}
The variable $\xi$ introduced in this lemma appears in the study of the
bias term. Indeed we have
\begin{eqnarray*}
  B_t & = & \e\left(\hat f_b(t)\right) -f(t) \\
  & = & \e\left(\frac{1}{n}\sum_{k=1}^n K_{t,b}(X_k)\right) \\
  & = & \int_0^1 K_{t,b}(x)f(x)dx - f(t) \\
  & = & \e\left(f(\xi)-f(t)\right).
\end{eqnarray*}
Now, we have to distinguish two cases depending on the position of
$\beta$ with respect to $1$.

\begin{enumerate}[{\bfseries a.}]
\item Assume that $0<\beta\leq 1$. Since $f$ belongs to $\Sigma(\beta,L)$,
  one can write:
  \begin{eqnarray*}
    |B_t| & \leq & \e\left(|f(\xi)-f(t)|\right) \\
    & \leq & L \e\left(|\xi-t|^\beta\right)\\
    & \leq & L \e\left(|\xi-t|\right)^\beta.
  \end{eqnarray*}
  Last inequality holds because $\beta\leq1$. Set $\eta>0$.
  \begin{eqnarray*}
    |B_t| & \leq & L \e\left(|\xi-t|\left(I_{|\xi-t|\leq\eta}+I_{|\xi-t|>\eta}\right)\right)^\beta\\
    & \leq & L \left(\eta+\sqrt{\sigma^2(\xi)\cdot \p\left(|\xi-t|>\eta \right)}\right)^\beta\\
    & \leq & L \left(\eta+\sqrt{\sigma^2(\xi)\cdot\frac{\sigma^2(\xi)}{\eta^2}}\right)^\beta\\
    & \leq & L \left(\eta+\frac{Cb}{\eta}\right)^\beta.
  \end{eqnarray*}
  Except for the last inequality which follows from
  Equation~(\ref{eq:21}) and~(\ref{eq:23}), all the others follow from
  classical probability inequalities. Optimizing in $\eta$, we obtain that $|B_t| \leq CLb^{\frac{\beta}{2}}$,
which allows to conclude.
\item Assume that $1<\beta\leq 2$. Using the mean-value theorem, we obtain:
  \begin{eqnarray*}
    |B_t| & = & \left|\e\left(f(\xi)-f(t)\right)\right| \\
    & = & \left|\e\left(f'(t)(\xi-t)+(\xi-t)(f'(\tilde\xi)-f'(t))\right)\right|,
  \end{eqnarray*}
  where $|\tilde\xi-t|\leq|\xi-t|$. Thus we obtain:
  \begin{eqnarray*}
    |B_t| & = &  \left|f'(t)(\e\left(\xi\right)-t)+\e\left((\xi-t)(f'(\tilde\xi)-f'(t))\right)\right|\\
    & \leq & |f'(t)|\cdot|\e\left(\xi\right)-t| + \e\left(|\xi-t|\cdot|f'(\tilde\xi)-f'(t)|\right) \\
    & \leq &  A(M+1)b + L\e\left(|\xi-t|^\beta\right),
  \end{eqnarray*}
  where $A=\sup_{f\in\Sigma(\beta,L)}\|f'\|_\infty<+\infty$ (see, for
  example, Tsybakov (2004), \S 1.2.1). Remark that last inequality follows from
  Equation~(\ref{eq:20}).

  Let us consider the term $\e\left(|\xi-t|^\beta\right)$. By
  introducing the mean of $\xi$ we obtain:
  \begin{eqnarray*}
    \e\left(|\xi-t|^\beta\right) & = & \e\left(|\xi-\e\left(\xi\right)+\e\left(\xi\right)-t|^\beta\right) \\
    & \leq & C \left(\e\left(|\xi-\e\left(\xi\right)|^\beta\right) +|\e\left(\xi\right)-t|^\beta\right).
  \end{eqnarray*}
  As $1<\beta\leq 2$, we obtain, using H\"older's inequality:
  \begin{eqnarray*}
    \e\left(|\xi-t|^\beta\right) & \leq & C \left(\sigma^\beta(\xi)+|\e\left(\xi\right)-t|^\beta\right)\\
    & \leq & C\left((Mb)^{\frac{\beta}{2}}+(Mb)^\beta\right)\\
    & \leq & Cb^{\frac{\beta}{2}}.
  \end{eqnarray*}
  Finally, we obtain again that $|B_t| \leq CLb^{\frac{\beta}{2}}$
that gives the lemma in the second case.
\end{enumerate}

\par 
{\bf 3.4 Proof of Lemma 2}
\par
Note that
$ Z_t  =  \hat f_b(t) - \e\left(\hat f_b(t)\right) =\frac{1}{n}\sum_{k=1}^n \eta_k$,
where the $\eta_k=K_{t,b}(X_k) - \e\left(K_{t,b}(X_k)\right)$'s are i.i.d. centered variables. Thanks to the following lemma (see
Bretagnolle and Hubert, 1979) it is possible to control precisely the
$p$-th moment of $Z_t$.
\begin{lemma}\label{lem:Bretagnolle}
  If $y_1,\ldots,y_n$ are $n$ i.i.d. variables such that
  $\e\left(y_1\right)=0$ and $\sigma^2(y_1)=v$, then the following
  inequalities hold:
  \begin{enumerate}
  \item If $p\leq 2$ then:
    \[
    \e\left(\left|\frac{1}{n}\sum_{k=1}^n y_k\right|^p\right) \leq
    \left(\frac{v}{n}\right)^{\frac{p}{2}}.
    \]
  \item If $p>2$ and moreover $\|y_1\|_\infty<+\infty$, then:
    \[
    \e\left(\left|\frac{1}{n}\sum_{k=1}^n y_k\right|^p\right) \leq
    C_p\left(\frac{v\|y_1\|_\infty^{p-2}}{n^{p-1}} +
      \left(\frac{v}{n}\right)^{\frac{p}{2}}\right).
    \]
  \end{enumerate}
\end{lemma}

This lemma will be apply with $y_k=\eta_k, (k=1,\ldots,n)$. Thus we
have to control two terms: $\sigma^2(\eta_1)$ on the one hand and
$\|\eta_1\|_{\infty}$ on the other hand.
\begin{enumerate}[{\bfseries a.}]
\item In order to control of $\sigma^2(\eta_1)$, let us compute:
  \begin{eqnarray*}
    \sigma^2(\eta_1)
    & \leq & \e\left(\left(K_{t,b}(X_1)\right)^2\right) \\
    & = & \int_0^1 \frac{x^{2\tb}(1-x)^{2\utb}}{B^2\!\left(\tb+1,\utb+1\right)}f(x)dx \\
    & = & A_b(t) \e[f(Z)],
  \end{eqnarray*}
  where:
  \begin{equation}\label{eq:A:t:b}
    A_b(t) = \frac{B\left(2\tb+1,2\utb+1\right)}{B^2\!\left(\tb+1,\utb+1\right)}
  \end{equation}
  and $Z$ has a beta distribution with parameters $2t/b+1$ and $2(1-t)/b+1$.
  Since $f$ belongs to $\Sigma(\beta,L)$ and thanks to the fact that
  the supremum of $\|f\|_{\infty}$ over $\Sigma(\beta,L)$ is finite,
  we obtain that
  $\sigma^2(\eta_1) \leq CA_b(t)$.
  Moreover, it is known (see Chen, 2000) that for $b$ small enough:
  \[
  A_b(t) \leq \frac{Cb^{-\frac{1}{2}}}{\sqrt{t(1-t)}}.
  \]
  Thus, for all $b$ small enough:
  \[
  \sigma^2(\eta_1) \leq \frac{Cb^{-\frac{1}{2}}}{\sqrt{t(1-t)}}.
  \]

\item Control of $\|\eta_1\|_{\infty}$. As $t$ is the mode of
  $K_{t,b}$ and thanks to Stirling's formula, it can be shown (see
  Chen, 2000) that:
  \[
  \sup_{x\in[0,1]} K_{t,b}(x) \leq
  \frac{Cb^{-\frac{1}{2}}}{\sqrt{t(1-t)}}.
  \]
  Thus we obtain:
  \begin{eqnarray*}
    \|\eta_1\|_\infty & = & \|K_{t,b}(X_1)-\e\left(K_{t,b}(X_1)\right)\|_\infty \\
    & \leq & \|K_{t,b}(X_1)\|_\infty + \left|\e\left(f(\xi)\right)\right|\\
    & \leq & \sup_{x\in[0,1]} K_{t,b}(x) + |f(t)| + |\e\left(f(\xi)\right)-f(t)| \\
    & \leq & \frac{Cb^{-\frac{1}{2}}}{\sqrt{t(1-t)}}+ \sup_{f\in\Sigma(\beta,L)}\|f\|_{\infty} +
    CLb^{\frac{\beta}{2}}.
  \end{eqnarray*}
  Thus, for small $b$ we have:
  \[
  \|\eta_1\|_\infty \leq \frac{Cb^{-\frac{1}{2}}}{\sqrt{t(1-t)}}.
  \]
\end{enumerate}
Now let us complete the proof of this lemma. We have to distinguish
two cases. The first one concerns the case where $p\leq2$. The second
one, the case where $p>2$.

First, let us assume that $p\leq2$. Applying
Lemma~\ref{lem:Bretagnolle} and using the bound on $\sigma^2(\eta_1)$
just obtained, we have:
\begin{eqnarray*}
  \e\left(|Z_t|^p\right) & = & \e\left(\left|\frac{1}{n}\sum_{k=1}^n
      \eta_k\right|^p\right)\\
  & \leq & \left( \frac{\sigma^2(\eta_1)}{n} \right)^{\frac{p}{2}} \\
   & \leq & C \left(\frac{1}{nb^{\frac{1}{2}}}\right)^{\frac{p}{2}} \left(t(1-t)\right)^{-\frac{p}{4}}.
\end{eqnarray*}

Last, if $p>2$, we have, thanks to Lemma~\ref{lem:Bretagnolle}:
\begin{eqnarray*}
  \e\left(|Z_t|^p\right) & = & \e\left(\left|\frac{1}{n}\sum_{k=1}^n \eta_k\right|^p\right)\\
  & \leq & C_p \left\{\frac{\sigma^2(\eta_1)\|\eta_1\|^{p-2}_\infty}{n^{p-1}}+\left(\frac{\sigma^2(\eta_1)}{n}\right)^{\frac{p}{2}}\right\}.
\end{eqnarray*}
Thanks to our bound on $\|\eta_1\|_\infty$ we obtain:
\begin{eqnarray*}
  \e\left(|Z_t|^p\right) & \leq & C_p \left\{\frac{\sigma^2(\eta_1)\sigma^{2p-4}(\eta_1)}{n^{p-1}}+\left(\frac{\sigma^2(\eta_1)}{n}\right)^{\frac{p}{2}}\right\} \\
  & \leq & C_p \left\{\left(\frac{\sigma^{2}(\eta_1)}{n}\right)^{p-1}+\left(\frac{\sigma^2(\eta_1)}{n}\right)^{\frac{p}{2}}\right\} \\
  & \leq & C_p \left\{\left(\frac{\sigma^{2}(\eta_1)}{n}\right)^{p-1-p/2}+1\right\}\left(\frac{\sigma^2(\eta_1)}{n}\right)^{\frac{p}{2}} \\
  & \leq & C \left(\frac{\sigma^2(\eta_1)}{n}\right)^{\frac{p}{2}}.
\end{eqnarray*}
Using our bound on $\sigma^2(\eta_1)$ we obtain:
\[
\e\left(|Z_t|^p\right) \leq
C \left(\frac{1}{nb^{\frac{1}{2}}}\right)^{\frac{p}{2}}
\left(t(1-t)\right)^{-\frac{p}{4}}.
\]
Taking all together lemma follows. Indeed, for any $p<4$ and $b$ small
enough we have:
\begin{equation*}
  \left(\int_0^1\e\left(|Z_t|^p\right)dt\right)^{\frac{1}{p}} \leq C \left(\frac{1}{nb^{\frac{1}{2}}}\right)^{\frac{1}{2}}
  \left(\int_0^1\left(t(1-t)\right)^{-\frac{p}{4}}dt\right)^{\frac{1}{p}}.
\end{equation*}

\par 
{\bf 3.5 Proof of Lemma~\ref{LEM:10}}

Let us consider a preliminary lemma which will be proved in Subsection~3.6.
\begin{lemma}\label{LEM:7}
For all $p\geq 2$ and $b\in(0,1)$ and all density $f$, we have:
\begin{equation}
\e\left(\|\hat f_b-f\|_p^{p}\right) \geq \int_0^1|B_t|^pdt\label{eq:minbiais}
\end{equation}
and:
\begin{equation}\label{eq:minvar}
\e\left(\|\hat f_b-f\|_p^{p}\right) \geq 2^{-p} \int_0^1\left(\e\left(Z_t^2\right)\right)^{\frac{p}{2}}dt.
\end{equation}
\end{lemma}

\par
\textit{Proof of i).}
In order to prove~(\ref{eq:minvar_0}) it is enough, thanks to~(\ref{eq:minvar}), to lower bound
$\int_0^1\left(\e_0\left(Z_t^2\right)\right)^{\frac{p}{2}}dt$.

Note that $\e_0\left(Z_t^2\right)$ can be written in the following way:
\begin{eqnarray*}
	\e_0\left(Z_t^2\right) & = & \frac{1}{n}\left(\e_0\left(\zeta_t^2\right)-\left(\e_0\left(\zeta_t\right)\right)^2\right)\\
	& = & \frac{1}{n}\left(\e_0\left(\zeta_t^2\right)-1\right)
\end{eqnarray*}
where $\zeta_t = K_{t,b}(X_1)$. Thus our goal is to minorate $\e_0\left(\zeta_t^2\right)$.
%
%
We have:
\begin{equation*}
	\e_0\left(\zeta_t^2\right) = \e\left(K_{t,b}^2(X_1)\right) = \int_0^1 K_{t,b}^2(x) dx =A_b(t)
\end{equation*}
where $A_b(t)$ is defined by Equation~(\ref{eq:A:t:b}).
%

Let us introduce the $R$-function defined as follows, for $z\geq 0$:
\[
	R(z) = \frac{1}{\Gamma(z+1)}\left(\frac{z}{e}\right)^{z}\sqrt{2\pi z}.
\]
It is well-known that $R$ is an increasing function such that $R(z)<1$ and, $R(z)\to 1$ as $z\to+\infty$.

Following Chen (1999) let us write $A_b(t)$ in terms of $R$-functions:
\[
		A_b(t) = c(b) \frac{1}{\sqrt{t(1-t)}} \frac{R^2(\frac{t}{b})R^2(\frac{1-t}{b})R(\frac{2}{t}+1)}{R(\frac{2t}{b})R(\frac{2(1-t)}{b})R^2(\frac{1}{b}+1)}
\]
where
$$c(b) = \frac{e^{-1}}{2\sqrt{\pi}}
		b\left(1+\frac{1}{b}\right)^{\frac{3}{2}}
		\left(1+\left(\frac{2}{b}+\frac{5}{2}\right)^{-1}\right)^{\frac{2}{b}+3}
		\geq  C b^{-\frac{1}{2}},$$
for $b$ small enough. Moreover, as $R$ is increasing we have:
\[
	\frac{R^2(\frac{t}{b})R^2(\frac{1-t}{b})R(\frac{2}{b}+1)}{R(\frac{2t}{b})R(\frac{2(1-t)}{b})R^2(\frac{1}{b}+1)} \geq \frac{R^2(\frac{t}{b})R^2(\frac{1-t}{b})}{R(\frac{2t}{b})R(\frac{2(1-t)}{b})R(\frac{1}{b}+1)}.
\]
Using the fact that $R(z)<1$, we obtain:
\[
	\frac{R^2(\frac{t}{b})R^2(\frac{1-t}{b})R(\frac{2}{b}+1)}{R(\frac{2t}{b})R(\frac{2(1-t)}{b})R^2(\frac{1}{b}+1)} \geq R^2\left(\frac{t}{b}\right)R^2\left(\frac{1-t}{b}\right).
\]

Thus we have:
\begin{eqnarray*}
	A_b(t) & \geq & C\frac{b^{-\frac{1}{2}}}{\sqrt{t(1-t)}} R^2\left(\frac{t}{b}\right)R^2\left(\frac{1-t}{b}\right) \\
	& \geq & C\frac{b^{-\frac{1}{2}}}{\sqrt{t(1-t)}} R^2\left(\frac{t}{b}\right)R^2\left(\frac{1-t}{b}\right) I_{[1/2, 1-b]}(t) \\
	& \geq & C\frac{b^{-\frac{1}{2}}}{\sqrt{t(1-t)}} R^2\left(\frac{1}{2b}\right)R^2(1) I_{[1/2, 1-b]}(t).
\end{eqnarray*}

For $b$ small enough, $R^2\left(\frac{1}{2b}\right)$ is greater than $R^2(1)$. Hence we obtain:
\begin{equation}\label{minoration_finale_Atb}
	A_b(t) \geq C\frac{b^{-\frac{1}{2}}}{\sqrt{t(1-t)}}.
\end{equation}

From Equation (\ref{minoration_finale_Atb}) we obtain:
\begin{eqnarray*}
	\e_0\left(\zeta_t^2\right) & \geq & Cb^{-\frac{1}{2}}\frac{1}{\sqrt{t(1-t)}}I_{[1/2, 1-b]}(t) \\
		& \geq & C b^{-\frac{1}{2}}\frac{1}{\sqrt{1-t}} I_{[1/2, 1-b]}(t).
\end{eqnarray*}

Thus, for $b$ small enough we obtain:
\begin{equation*}
	\e_0\left(Z_t^2\right) \geq C \frac{1}{nb^{1/2}}\frac{1}{\sqrt{1-t}}I_{[1/2, 1-b]}(t) -\frac{1}{n}.
\end{equation*}

Finally, we can write:
\begin{eqnarray*}
	\int_0^1\left(\e_0\left(Z_t^2\right)\right)^{\frac{p}{2}} dt & \geq & C \left(\frac{1}{nb^{1/2}}\right)^{\frac{p}{2}} d_n(b,p). \\
\end{eqnarray*}

\par
\textit{Proof of ii).}
We define for $k=1\ldots,N$, $t_k = (2k-1)/{(4N)}$ and the following intervals: $T_k =[t_k-\varepsilon b^{1/2}, t_k+\varepsilon b^{1/2}]$, $I_k = [t_k-b^{1/2}, t_k+b^{1/2}]$ and $J_k=[t_k-1/(4N), t_k+1/(4N))$ where $0<\varepsilon<1/2$ will be chosen. Let us recall that $1/(4N)$ is very close to $5b^{1/2}$. Note that we have the following inclusions: $T_k\subset I_k \subset J_k$.

Now, let us lower bound the integrated bias:
\begin{eqnarray*}
	\int_0^1 \left|B_t\right| dt & \geq & \sum_{\ell=1}^N \int_{T_{2\ell-1}}\left|\int_0^1 K_{t,b}(x)\left(f_\beta(t)-f_\beta(x)\right) dx\right| dt\\
	& \geq & \sum_{\ell=1}^N \int_{T_{2\ell-1}} (A_\ell(t)-B(t)) dt
\end{eqnarray*}
where
\[
	A_{\ell}(t) = \int_{J_{2\ell-1}\cap I_{2\ell-1}^{c}}K_{t,b}(x)\left(f_\beta(t)-f_\beta(x)\right) dx
\]
and
\[
	B(t) = \int_{\{f_\beta(x)\geq f_\beta(t)\}}K_{t,b}(x)\left(f_\beta(x)-f_\beta(t)\right) dx.
\]

Since for $t\in T_{2\ell-1}$ and $x\in J_{2\ell-1}\cap I_{2\ell-1}^{c}$ we have:
\begin{eqnarray*}
	f_\beta(t)-f_\beta(x) & = & L_\beta\left(\left|x-t_\ell\right|^\beta-\left|t-t_\ell\right|^\beta\right) \\
		& \geq & L_\beta\left(b^{\beta/2}-\varepsilon^{\beta}b^{\beta/2}\right) \\
		& \geq & C b^{\beta/2},
\end{eqnarray*}
it follows that:
\begin{eqnarray*}
	A_{\ell}(t) & \geq & C b^{\beta/2}\int_{J_\ell\cap I_\ell^{c}}K_{t,b}(x)dx\\
		& = & C b^{\beta/2}\p\left(\xi\in J_{2\ell-1}\cap I_{2\ell-1}^{c}\right),
\end{eqnarray*}
where $\xi\sim K_{t,b}$. This probability can be estimated:
\begin{eqnarray*}
	\p\left(\xi\in J_{2\ell-1}\cap I_{2\ell-1}^{c}\right)
	& \geq & \p\left(t_{2\ell-1}+b^{1/2}<\xi<t_{2\ell-1}+5b^{1/2}\right)\\
	& \geq & \p\left(t_{2\ell-1}-t+t+b^{1/2}<\xi<t_{2\ell-1}-t+t+5b^{1/2}\right)\\
	& \geq & \p\left(t+2b^{1/2}<\xi<t+4b^{1/2}\right)\\
	& \geq & 2b^{1/2} K_{t,b}\left(t+4b^{1/2}\right)
\end{eqnarray*}

For $1/4\leq t\leq 3/4$ we have:
\[
	K_{t,b}(t+4b^{1/2}) \mathop{\sim}_{b\to0} K_{t,b}(t) \exp\left(-\frac{8}{t(1-t)}\right)\geq Cb^{-1/2}.
\]
Indeed, using the $R$ function we obtain:
\begin{eqnarray*}
	K_{t,b}(t) & = & \frac{R\left(\tb\right)R\left(\utb\right)\sqrt{b}(1+1/b)}{R\left(1/b\right)\sqrt{2\pi t(1-t)}}\\
	 & \geq & C b^{-1/2}
\end{eqnarray*}
for $1/4\leq t\leq 3/4$ and $b$ small enough.
Thus, for $\ell$ such that $T_{2\ell-1}\subset [1/4,3/4]$, $t\in T_{2\ell-1}$ and for $b$ small enough, we have:
\[
	A_\ell(t) \geq Cb^{\beta/2},
\]
where $C$ is an absolute constant independent of $b$ and $\varepsilon$.

For $t\in T_{2\ell-1}$ we have:
\begin{eqnarray*}
	B(t) & = & \int_{\bigcup_{k=1}^N\{|x-t_{2k-1}|\leq |t-t_{2\ell-1}|\}}K_{t,b}(x)\left(f_\beta(x)-f_\beta(t)\right) dx\\
		& \leq & \sum_{k=1}^N \int_{\{|x-t_{2k-1}|\leq |t-t_{2\ell-1}|\}}K_{t,b}(x) L_\beta(|t-t_{2\ell-1}|^\beta-|x-t_{2k-1}|^{\beta}) dx\\
		& \leq & C\varepsilon^\beta b^{\beta/2}.
\end{eqnarray*}

Taking all together we obtain the following lower bound:
\begin{eqnarray*}
	\int_0^1 |B_t| dt & \geq & \mathop{\sum_{\ell=1}^N}_{T_{2\ell-1}\subset [1/4,3/4]} \int_{T_{2\ell-1}} A_\ell(t) dt - \sum_{\ell=1}^N \int_{T_{2\ell-1}}B(t) dt\\
	& \geq & CN\varepsilon b^{1/2} b^{\beta/2} - N \varepsilon^{\beta+1} b^{1/2} b^{\beta/2} \\
	& \geq & b^{\beta/2}\left(c_1\varepsilon-c_2\varepsilon^{\beta+1}\right)
\end{eqnarray*}
where $c_1$ and $c_2$ are absolute positive constants. Thus, for $\varepsilon = \min(1/2,(c_1/c_2)^{1/\beta}/2)$ we obtain:
\[
	\int_0^1 |B_t| dt \geq C b^{\beta/2},
\]
and thus, using the convexity of $x\mapsto x^p$,
\[
	\int_0^1 |B_t|^p dt \geq C b^{\beta p/2}.
\]
Using~(\ref{eq:minbiais}), this achieves the proof of ii).

\par

\par
\textit{Proof of iii).}
For $b$ small enough, we have:
\begin{eqnarray*}
	B_t & = & \e_3\left(K_{t,b}(X_1)\right) - f_3(t) \\
		& = & \e\left(f_3(\xi)\right) -f_3(t) \\
		& = & 2(\e\left(\xi\right)-t) \\
		& \geq & Cb(1-2t),
\end{eqnarray*}
where $\xi$ is defined as in Lemma~\ref{LEM:JKB} and last inequality is a consequence of (\ref{eq:20}). As $|1-2t|^p$ is integrable on $[0,1]$, it follows:
\[
	\int_0^1 |B_t|^pdt \geq C b^p.
\]
Using Equation~(\ref{eq:minbiais}) we deduce Equation~(\ref{eq:minbiais_3}).

\par 
{\bf 3.6 Proof of Lemma~\ref{LEM:7}}

Since $x\mapsto|x|^p$ is a convex function
$$\left|B_t\right|^p=\left|\e\left(\hat f_b(t)-f(t)\right)\right|^p\le \e\left(\left|\hat f_b(t)-f(t)\right|^p\right),$$
which implies~(\ref{eq:minbiais}) taking the integral on $[0,1]$ in both sides of the last inequality.

Since $\|Z_t\|_p\le \|B_t\|_p+\|\hat f_b-f\|_p$, we deduce, using~(\ref{eq:minbiais}) that, for all function~$f$,
$$\e\left(\|Z_t\|_p^p\right)\le 2^{p-1}\e\left(\|B_t\|_p^p+\|\hat f_b-f\|_p^p\right)\le 2^p\e\left(\|\hat f_b-f\|_p\right).$$
Then, we have
\begin{eqnarray*}
	\e\left(\|\hat f_b-f\|_p^{p}\right) & \geq & 2^{-p} \e\left(\|Z_t\|_p^p\right) \\
		& = & 2^{-p} \e\left(\int_0^1|Z_t|^p dt\right) \\
		& = & 2^{-p}\int_0^1\e\left(|Z_t|^p\right)dt \\
		& \geq & 2^{-p} \int_0^1\left(\e\left(Z_t^2\right)\right)^{\frac{p}{2}}dt.
\end{eqnarray*}
This proves~(\ref{eq:minvar}).

\par
\par

\noindent {\large\bf Acknowledgment}
\par
Karine Bertin and Nicolas Klutchnitkoff are supported by Project Fondecyt 1090285. Karine Bertin has been supported by project Laboratory ANESTOC PBCT ACT 13.
\par

\par

\noindent{\large\bf References}
\begin{description}
\item Bretagnolle, J.and Huber, C. (1979). Estimation des densités: risque minimax. {\it Z. Wahrsch. Verw. Gebiete},  {\bf 47}, no. 2, 119-137.
\item Bouezmarni, T. and van Bellegem (2009). Nonparametric beta kernel estimator for long memory time series, {\it technical report}.
\item Bouezmarni, T. and Rolin, J-M. (2003). Consistency of the beta kernel density function estimator,
{\it Canad. J. Statist.}, {\bf 31} , no. 1, 89-98.
\item Charpentier, A. and  Oulidi A. (2010). Beta kernel quantile estimators of heavy-tailed loss distributions. {\it Statistics and computing}, {\bf 20}, no. 1, 35-55.
\item
Chen, S. X. (1999). A beta kernel estimator for density functions with compact supports. {\it Comput. Statist. Data Anal.}, {\bf 31}, 131-145.
\item
Chen, S. X. (2000). Beta kernel smoothers for regression curves. {\it Statistica Sinica}, {\bf 10}, 73-91.
\item
Ibragimov, I. A.; Hasminskii, R. Z. (1981). {\it Statistical estimation. Asymptotic theory.} Springer-Verlag, New York-Berlin.
\item
Johnson, N. L., Kotz, S. and Balakrishnan, N. (1994). {\it Continuous Univariate Distributions.} Wiley, New York.
\item
Tsybakov, A. (2004). {\it Introduction \`a l'estimation non param\'etrique.}  Springer, Berlin.
\end{description}

\vskip .65cm
\noindent
Karine Bertin, Departamento Estad\'{i}stica, CIMFAV, Universidad de Valpara\'{i}so, Avenida Gran Breta\~na 1091, Playa Ancha, Valpara\'{i}so, Chile, tel/fax: 0056322508268
\vskip 2pt
\noindent
E-mail: karine.bertin@uv.cl
\vskip 2pt
\noindent
Nicolas Klutchnitkoff,
Institut de Recherche Mathématique Avancée (IRMA),
CNRS : UMR7501 – Université de Strasbourg,
7 rue René descartes, Srasbourg,
France, tel: 0033368850186
\vskip 2pt
\noindent
E-mail: klutchni@math.unistra.fr
\vskip .3cm
\end{document}